\documentclass[preprint,12pt]{article}
\usepackage{amssymb}
\usepackage{mathrsfs}
\usepackage{amsfonts}
\usepackage{graphicx}
\usepackage{colortbl,dcolumn}
\usepackage{amsmath}
\usepackage{psfrag}
\usepackage{booktabs}
\usepackage{slashbox}
\usepackage{cite}
\numberwithin{equation}{section}
\usepackage[top=1.2in, bottom=1.2in, left=1in, right=1in]{geometry}

\newtheorem{theorem}{Theorem}[section]

\begin{document}
\title{Modified equations for weak stochastic symplectic schemes via their generating functions }
       \author{Lijin Wang\footnotemark[1],
        Jialin Hong\footnotemark[2],\\
        {\small \footnotemark[1] School of Mathematical Sciences, University of Chinese Academy of Sciences,}\\
    {\small Beijing 100049, P.R.China }\\
       {\small \footnotemark[2] Institute of Computational Mathematics and Scientific/Engineering Computing,}\\
       {\small Academy of Mathematics and Systems Science, Chinese Academy of Sciences, }\\
         {\small Beijing 100190, P.R.China }} 
\maketitle

       \begin{abstract}
          {\rm\small In this paper, a systematic approach of constructing modified equations for weak stochastic symplectic methods of stochastic Hamiltonian systems is given via using the generating functions of the stochastic symplectic methods. This approach is valid for stochastic Hamiltonian systems with either additive noises or half-multiplicative noises, and we prove that the modified equation of the weak stochastic symplectic methods are perturbed stochastic Hamiltonian systems of the original systems, which reveals in certain sense the reason for the good long time numerical behavior of stochastic symplectic methods.}\\


\end{abstract}

\section{Introduction}
\label{sect:intro}

\label{sec;introduction}
\label{intro}
The construction of modified equation is generally the first step of the backward error analysis for numerical solution of differential equations, which has obtained much success during the last decades (\cite{ei,feng1,fied,gri}), especially in the qualitative analysis of symplectic methods for Hamiltonian systems (\cite{bene,hairer1,moser,mu,reich,sanz1,sanz,tang,yo}). It constructs perturbed Hamiltonian systems to approximate the symplectic methods, reveals the reason for their long-time good behavior, and implies construction of higher order methods. For a more detailed review of these see e.g. the monographs \cite{hairer,leim}.

The extension of modified equation to stochastic context is a relatively new topic. To knowledge of the authors, some prior work include \cite{ab,fa,18,sh,wwang,zy}, and so on. In \cite{wwang}, linear Langevin equations are considered. In \cite{sh}, modified equations that approximate the forward and backward Euler methods in the sense of weak convergence for It\^{o} stochastic differential equations (SDEs) with additive noise are given, with discussion of extension to other type of equations and approximation. \cite{zy} describes a general framework for deriving modified equations for SDEs with respect to weak convergence, using weak stochastic Taylor expansions. \cite{fa} builds modified equations not at the level of the SDEs, but at the level of the Kolmogorov generators associated with the process solution of the SDEs.
\cite{18} uses modified equations to exhibit the poor behavior of the Euler methods for small random perturbations of Hamiltonian flows, and \cite{ab} proposes a new method for constructing numerical integrators with high weak order for the time integration of SDEs, inspired by the theory of stochastic modified equations.

In this paper, we attempt to construct modified equations for weak symplectic methods for stochastic Hamiltonian systems in terms of their generating functions, which is an extension of the approach dealing with symplectic methods for deterministic Hamiltonian systems \cite{bene,mu}. We prove that the modified equation of weak symplectic methods are again perturbed stochastic Hamiltonian systems, and describe a systematic approach of establishing modified equations that approximate the symplectic methods to certain weak order, using the stochastic generating functions of these symplectic schemes.

The contents are arranged as follows. Section 2 is an introduction of the stochastic modified equations theory. Section 3 reviews the stochastic generating function theory for stochastic symplectic methods. In section 4 we develop the approach of constructing modified equations of weak symplectic methods via their generating functions, and prove that the modified equations of weak symplectic methods are perturbed stochastic Hamiltonian systems of the original systems. Section 5 gives some examples, and Section 6 is a brief conclusion.

\section{Stochastic modified equations}
\label{sec:1}
Given a stochastic differential equation
\begin{equation}\label{2.1}
  dX(t)=a(X)dt+\sum_{r=1}^m\sigma_r(X)d W_r(t),\quad X(0)=x,
\end{equation}
where the drift and diffusion coefficients $a(X)$ and $\sigma_r(X)$ are $\mathbf{R}^d\rightarrow\mathbf{R}^d$ functions satisfying conditions that guarantee the existence and uniqueness of the solution (see e.g. \cite{mao,oksen} for details). $W(t)=(W_1(t),\cdots, W_m(t))^T$ is $m$-dimensional standard Wiener process.  A numerical approximation of it $X_0,X_1,\cdots$ with step size $h$ is said to converge weakly of order $p$, if for any $T>0$ with $N_T h=T$
\begin{equation}\label{2.2}
 |\mathbf{E}(\phi(X_{N_T}))-\mathbf{E}(\phi(X(T)))|=O(h^p)
\end{equation}
for $\phi\in C_P^{2p+2}(\mathbf{R}^d,\mathbf{R})$, the space of $2p+2$ times continuously differentiable functions mapping from $\mathbf{R}^d$ to $\mathbf{R}$ which, together with their partial derivatives up to order $2p+2$ have polynomial growth. As test functions, it is usually enough to consider polynomials up to order $2p+1$.

In the aforementioned literature, the $q$-th order stochastic modified equation of a weak numerical method is the modified SDE
\begin{equation}\label{2.3}
d\tilde{X}=A(\tilde{X})dt+\sum_{r=1}^m \Gamma_r(\tilde{X})dW_r(t)
\end{equation}
with
\begin{equation}\label{2.5}
\begin{array}{l}
  A(\tilde{X})=a(\tilde{X})+a_1(\tilde{X})h+a_2(\tilde{X})h^2+\cdots+a_q(\tilde{X})h^q\\
  \Gamma_r(\tilde{X})=\sigma_r(\tilde{X})+\sigma_{r,1}(\tilde{X})h+\cdots+\sigma_{r,q}(\tilde{X})h^q,
\end{array}
\end{equation}
where the functions $a_i$ and $\sigma_{r,i}$ are to be determined such that for any $T>0$ with $T=hN_T$
\begin{equation}\label{2.4}
|\mathbf{E}(\phi(X_{N_T}))-\mathbf{E}(\phi(\tilde{X}(T)))|=O(h^{p+q'}),
\end{equation}
for some $q'>0$, meaning that the solution of (\ref{2.3}) is weakly $q'$-order 'closer' to the numerical solution than the SDE (\ref{2.1}). In the ODE case, $q'$ can be infinity for any $T=hN_T$, $N_T\in \mathbf{Z}_{+}\bigcup\{0\}$, meaning that the discrete numerical solution locates exactly on the continuous true solution curve of its modified equation.

In \cite{zy}, the functions $a_i$ and $\sigma_{r,i}$ are found by using the Taylor expansion of $\mathbf{E}(\phi(X)|X(0)=x)$, i.e., the Taylor expansion of the solution to the backward Kolmogorov equation associated with the SDE (\ref{2.1})
\begin{equation}\label{2.6}
\frac{\partial u}{\partial t}=\mathcal{L}_0 u,\quad u(x,0)=\phi(x),
\end{equation}
where, in the case of $d=m=1$ in (\ref{2.1}), $\mathcal{L_0}u:=a(x)\frac{du}{dx}+\frac{1}{2}\sigma^2(x)\frac{d^2 u}{dx^2}$. The solution of probabilistic
sense of (\ref{2.6}) is just (\cite{oksen})
\begin{equation}\label{2.7}
u(x,t)=\mathbf{E}(\phi(x(t))|x(0)=x),
\end{equation}
for which it holds
\begin{equation}\label{2.8}
u(x,h)-\phi(x)=\sum_{k=0}^N \frac{h^{k+1}}{(k+1)!}\mathcal{L}_0^{k+1}\phi(x)+O(h^{N+2})
\end{equation}
if $u$ is $N+1$ times differentiable with respect to $t$.

The expectation of $\phi$ of a one-step numerical approximation of weak order $p$, $\mathbf{E}(\phi(x(h))|x(0)=x)=:u_{num}(x,h)$ should have expansion coinciding with (\ref{2.8}) up to terms of order $O(h^p)$ which corresponds a local error of order $p+1$. On the other hand, there is the backward Kolmogorov equation associated with the modified equation (\ref{2.3})
\begin{equation}\label{2.9}
\frac{\partial u_{mod}}{\partial t}=\mathcal{L}^h u_{mod},\quad u_{mod}(x,0)=\phi(x),
\end{equation}
with $\mathcal{L}^h u=A(x)\frac{du}{dx}+\frac{1}{2}\Gamma^2(x)\frac{d^2u}{dx^2}$ for $d=m=1$. To obtain a modified equation of order $q$, one needs to equate the expansion of $u_{num}(x,h)$ which is known from the numerical method itself, and the expansion of $u_{mod}(x,h)$ which contains unknown functions to be determined, up to terms of order $O(h^{p+q})$. In this way the assumed unknown functions in the modified equations can be found (\cite{zy}).

In this paper, we construct modified equations for weak symplectic methods, not based on the backward Kolmogorov equations associated with the stochastic Hamiltonian systems, but via employing the stochastic generating functions that produce the weak symplectic schemes which are introduced in the following section.
\section{Generating functions for weak stochastic symplectic methods}
\label{sec:2}
Given a stochastic Hamiltonian system
\begin{equation}\label{3.1}
\begin{array}{l}
dp=-H_qdt-\sum_{r=1}^m H_{r_{q}}\circ dW_r(t),\quad p(0)=p_0,\\
dq=H_pdt+\sum_{r=1}^m H_{r_{p}}\circ dW_r(t),\quad q(0)=q_0,
\end{array}
\end{equation}
where $H(p,q)$, $H_r(p,q)$ are Hamiltonian functions, $(W_1(t),\cdots, W_m(t))^T$ is $m$-dimensional standard Wiener process, and $H_q$, $H_p$ denotes the partial derivatives of $H$ with respect to $q$ and $p$, respectively. The same holds for $H_r$. The small circle '$\circ$' before $dW_r(t)$
denotes SDEs of Stratonovich sense.

It is known that this system possesses the symplectic structure (\cite{mil1})
\begin{equation}\label{3.2}
dp(t)\wedge dq(t)=dp_0\wedge dq_0,\quad \forall t\geq 0.
\end{equation}
A numerical discretization $(p_n,q_n)_n$ that preserves this structure is called a symplectic method, characterized by
\begin{equation}\label{3.3}
dp_{n+1}\wedge dq_{n+1}=dp_n\wedge dq_n,\quad \forall n\in \mathbf{Z},\,\,\,n\geq 0.
\end{equation}

Symplectic methods for stochastic Hamiltonian systems have been shown to possess much superiorities than non-symplectic ones in long time simulation of stochastic Hamiltonian systems (\cite{milbook,mil2}). To construct such methods, the stochastic generating function approach was established (\cite{anton,deng,wang,dcds}). It is based on the fact that each symplectic mapping $(p_n,q_n)\mapsto (p_{n+1},q_{n+1})$ can be associated with a generating function (\cite{fengg,hairer,deng,wang}), e.g., $S^1(p_{n+1},q_n)$, such that
\begin{equation}\label{3.4}
\begin{array}{l}
p_{n+1}=p_n-\frac{\partial S^1}{\partial q_n},\\
q_{n+1}=q_n+\frac{\partial S^1}{\partial p_{n+1}}.
\end{array}
\end{equation}

There are other generating functions such as $S^2(p_n, q_{n+1})$ and $S^3(\frac{p_n+p_{n+1}}{2},\frac{q_n+q_{n+1}}{2})$, as illustration we only consider $S^1$ in this article. Methods based on the other generating functions are similar.

It can be proved that (\cite{deng,wang,dcds}), the generating function $S^1(P,q,t)$ which generates almost surely the phase flow of the stochastic Hamiltonian system (\ref{3.1}) satisfies the following stochastic Hamilton-Jacobi partial differential equation
\begin{equation}\label{3.5}
d_t S^1=H(P,q+\frac{\partial S^1}{\partial P})dt+\sum_{r=1}^m H_r(P,q+\frac{\partial S^1}{\partial P})\circ dW_r(t),\quad S^1(P,q,0)=0.
\end{equation}
On the other hand, every sufficiently smooth $S^1(P,q,t)$ satisfying (\ref{3.5}) can generate the phase flow of the stochastic Hamiltonian system (\ref{3.1}) $\varphi_t: (p,q)\mapsto (P(t),Q(t))$ almost surely via the relation
\begin{equation}\label{3.6}
p=P(t)+\frac{\partial S^1}{\partial q}(P(t),q,t),\quad Q(t)=q+\frac{\partial S^1}{\partial P}(P(t),q,t).
\end{equation}

Assume that $S^1(P,q,t)$ has a series expansion
\begin{equation}\label{3.7}
S^1(P,q,t)=\sum_{\alpha}G_{\alpha}(P,q)J_{\alpha},
\end{equation}
where $\alpha=(j_1,j_2,\cdots,j_l)$ denotes the multi-index of stochastic multiple integrals, $j_i\in\{0,1,\cdots,m\}$ $(i=1,\cdots,l)$, $l\geq 1$, and
\begin{equation}\label{3.8}
J_{\alpha}=\int_0^t\int_0^{s_l}\cdots\int_0^{s_2}\circ dW_{j_1}(s_1)\circ dW_{j_2}(s_2)\circ\cdots\circ dW_{j_l}(s_l).
\end{equation}
For convenience denote $ds=dW_{0}(s)$.

Substituting the series expansion (\ref{3.7}) into the stochastic Hamilton-Jacobi PDE (\ref{3.5}), performing Taylor expansion of $H$ and $H_r$ $(r=1,\cdots,m)$ at $(P,q)$ in (\ref{3.5}), and then equating coefficients of like powers on both sides of the equation, one obtains (\cite{anton,deng})
\begin{equation}\label{3.9}
G_{\alpha}=\sum_{i=1}^{l(\alpha)-1}\frac{1}{i!}\sum_{k_1,\cdots,k_i=1}^d\frac{\partial^i H_r}{\partial q_{k_1}\cdots \partial q_{k_i}}\sum_{\begin{array}{c}l(\alpha_1)+\cdots+l(\alpha_i)=l(\alpha)-1\\ \alpha-\in\Lambda \alpha_1,\cdots,\alpha_i\end{array}}\frac{\partial G_{\alpha_1}}{\partial p_{k_1}}\cdots\frac{\partial G_{\alpha_i}}{\partial p_{k_i}}
\end{equation}
for $\alpha=(i_1,\cdots,i_{l-1},r)$ with $l>1$, $i_1,\cdots,i_{l-1},r\in\{0,1,\cdots,m\}$, and there is no duplicate in $\alpha$. If there are duplicates in $\alpha$, one can still use the formula after assigning different subscripts to the duplicates. For $l(\alpha)=1$, i.e., $\alpha=(r)$, then
\begin{equation}\label{3.9.5}
G_r=H_r.
\end{equation}
Note that, $l(\alpha)$ denotes the length of $\alpha$, and $\alpha-$ is the multi-index resulted from discarding the last index of $\alpha$. $\Lambda_{\alpha_1,\cdots, \alpha_k}$ can be defined recursively as follows (\cite{anton,deng}). First, define the concatenation $'*'$ of the indices $\alpha=(j_1,\cdots,j_l)$ and $\alpha'=(j_1',\cdots,j_{l'}')$ as $\alpha\ast \alpha'=(j_1,\cdots,j_l,j_1',\cdots,j_{l'}')$. The concatenation of a set of multi-indices $\Lambda$ and $\alpha$ is $\Lambda\ast\alpha=\{\beta*\alpha|\beta\in \Lambda\}$. Then, define
\begin{equation}\label{3.10}
\Lambda_{ \alpha_1, \alpha_2}=\left\{\begin{array}{l}\{(j_1,j_1'),(j_1',j_1)\},\,\,\,\, \mbox{if}\,\,\,\, l=l'=1\\ \{\Lambda _{(j_1),\alpha_2- }*(j_{l'}'),\alpha_2*(j_1)\}, \,\,\,\,\mbox{if}\,\,\,\, l=1, l'\neq 1\\ \{\Lambda_{ \alpha_1-,(j_1')}*(j_l),\alpha_1 *(j_{1}')\},\,\,\,\,\mbox{if}\,\,\,\,l\neq1, l'=1\\ \{\Lambda_{\alpha_1-,\alpha_2}*(j_l),\Lambda_{\alpha_1,\alpha_2-}*(j_{l'}')\} \,\,\,\,\mbox{if}\,\,\,\, l\neq 1, l'\neq 1\end{array}\right.
\end{equation}
For $k>2$, define $\Lambda_{\alpha_1,\cdots, \alpha_k}=\{\Lambda_{\beta, \alpha_k}|\beta\in\Lambda_{\alpha_1,\cdots, \alpha_{k-1}}\}$.

With the $G_{\alpha}$ given in (\ref{3.9}), replacing the $t$ by $h$, $P$ by $p_{n+1}$, and $q$ by $q_n$ in (\ref{3.7}), and truncating the series (\ref{3.7}) to certain terms, a symplectic scheme of corresponding order can be obtained via the relation (\ref{3.4}).

To obtain a symplectic scheme of weak order $k$, one can firstly transform the multiple Stratonovich integrals $J_{\alpha}$ in (\ref{3.7}) to Ito integrals
\begin{equation}\label{3.11}
I_{\alpha}:=\int_0^t\int_0^{s_l}\cdots\int_0^{s_2}dW_{j_1}(s_1) dW_{j_2}(s_2)\cdots dW_{j_l}(s_l),
\end{equation}
using the relations (\cite{kloeden})
\begin{equation}\label{3.12}
J_{\alpha}=I_{(j_l)}[J_{\alpha-}]+\chi_{\{j_l=j_{l-1}\neq0\}}I_{(0)}[\frac{1}{2}J_{(\alpha-)-}],\quad l(\alpha)\geq 2,
\end{equation}
where $\chi_A$ denotes the indicator function of set $A$, and for $l(\alpha)=1$, $J_{\alpha}=I_{\alpha}$. Then, one should include in the truncation of the series (\ref{3.7}) all terms with index $\alpha$ which satisfies $l(\alpha)\leq k$ (\cite{kloeden,anton,deng}).

In Section 5 some examples of weak symplectic schemes for stochastic Hamiltonian systems produced by the generating functions are illustrated.
\section{Modified equations for weak symplectic schemes via their generating functions}
Inspired by the modified equation (\ref{2.3}) with (\ref{2.5}), we prove the following theorem about the modified equation of weak symplectic methods for stochastic Hamiltonian systems (\ref{3.1}).
\begin{theorem}
Given the stochastic Hamiltonian system (\ref{3.1}), and a weak symplectic scheme $\psi_h:(p,q)\mapsto (P,Q)$ with generating function $S^1(P,q,h)=\sum_{\alpha\in \Lambda_{\psi}}G_{\alpha}(P,q)J_{\alpha}^h$ where $G_{\alpha}(P,q)$ are defined on an open set $D$, and $J_{\alpha}^h$ are multiple stochastic integrals as in (\ref{3.8}) but defined on $(0,h)$. $\Lambda_{\psi}$ is the set of indices associated with $\psi_h$. Then the modified equation of the weak symplectic scheme $\psi_h$ is a stochastic Hamiltonian system
\begin{equation}\label{4.1}
\begin{array}{l}
dp=-\tilde{H}_{0_{q}}dt-\sum_{r=1}^m \tilde{H}_{r_q}\circ dW_r(t), \quad p(0)=p_0,\\
dq=\tilde{H}_{0_{p}}dt+\sum_{r=1}^m \tilde{H}_{r_{p}}\circ dW_r(t),\quad q(0)=q_0,
\end{array}
\end{equation}
where
\begin{equation}\label{4.2}
\begin{array}{l}
\tilde{H}_0=H+H_0^{[1]}h+H_0^{[2]}h^2+\cdots,\\
\tilde{H}_r=H_r+H_r^{[1]}h+H_r^{[2]}h^2+\cdots,
\end{array}
\end{equation}
and $H_j^{[i]}(p,q)$ are defined on $D$ for $j=0,1,2,\cdots$, and $i=1,2,\cdots$.
\end{theorem}
$Proof.$\\ \\
For convenience, denote $H=H_0^{[0]}$, $H_r=H_r^{[0]}$. The proof is also the process of finding the unknown functions $H_j^{[i]}$, for $j=0,1,2,\cdots$, and $i=1,2,\cdots$. Suppose the generating function that generates the stochastic Hamiltonian system (\ref{4.1}) is $\tilde{S}^1(P,q,t)$ which satisfies the stochastic Hamilton-Jacobi PDE
\begin{equation}\label{4.3}
d\tilde{S}^1=\tilde{H}_0(P,q+\frac{\partial \tilde{S}^1}{\partial P})dt+\sum_{r=1}^m \tilde{H}_r(P,q+\frac{\partial \tilde{S}^1}{\partial P})\circ dW_r(t), \quad \tilde{S}^1(P,q,0)=0.
\end{equation}
Assume the solution of (\ref{4.3}) has the form
\begin{equation}\label{4.4}
\tilde{S}^1(P,q,t)=\sum_{\alpha}\tilde{G}_{\alpha}(P,q,h)J_{\alpha}.
\end{equation}
According to the formula (\ref{3.9}), we have for $\alpha=(i_1,\cdots,i_{l-1},r)$, with $l\geq 2$
\begin{equation}\label{4.5}
\tilde{G}_{\alpha}=\sum^{l(\alpha)-1}_{i=1}\frac{1}{i!}\sum_{k_1,\cdots,k_i=1}^d\frac{\partial^i \tilde{H}_r}{\partial q_{k_1}\cdots \partial q_{k_i}}
\sum_{\begin{array}{c} l(\alpha_1)+\cdots+l(\alpha_i)=l(\alpha)-1\\ \alpha-\in\Lambda \alpha_1,\cdots,\alpha_i\end{array}}\frac{\partial \tilde{G}_{\alpha_1}}{\partial p_{k_1}}\cdots\frac{\partial \tilde{G}_{\alpha_i}}{\partial p_{k_i}},
\end{equation}
and $\tilde{G}_r=\tilde{H}_r$ for $\alpha=(r)$. Now suppose
\begin{equation}\label{4.6}
\tilde{G}_{\alpha}(P,q,h)=G_{\alpha}^{[0]}(P,q)+G_{\alpha}^{[1]}(P,q)h+G_{\alpha}^{[2]}h^2+\cdots.
\end{equation}
Substituting the series of $\tilde{H}_r$ in (\ref{4.2}), and that of $\tilde{G}_{\alpha_j}$ $(j=1,\cdots,i)$ as in (\ref{4.6}) into the right hand side of (\ref{4.5}), using (\ref{4.6}) as the left hand side of (\ref{4.5}), and then compare like powers of $h$ on both sides of (\ref{4.5}), we obtain for $\alpha=(i_1,\cdots,i_{l-1},r)$ with $l\geq 2$
\begin{eqnarray}\label{4.7}
G_{\alpha}^{[k]}=&\sum^{l(\alpha)-1}_{i=1}\frac{1}{i!}\sum_{k_1,\cdots,k_i=1}^d\sum_{j+j_1+\cdots+j_i=k}\frac{\partial^i H_r^{[j]}}{\partial q_{k_1}\cdots \partial q_{k_i}}\cdot\nonumber\\
&\sum_{\small \begin{array}{c} l(\alpha_1)+\cdots+l(\alpha_i)=l(\alpha)-1\\ \alpha-\in\Lambda \alpha_1,\cdots,\alpha_i\end{array}}\frac{\partial G_{\alpha_1}^{[j_1]}}{\partial p_{k_1}}\cdots\frac{\partial G_{\alpha_i}^{[j_i]}}{\partial p_{k_i}},
\end{eqnarray}
and for $\alpha=(r)$,
\begin{equation}\label{4.7.5}
G_r^{[k]}=H_r^{[k]}.
\end{equation}

According to (\ref{4.4}) and (\ref{4.6}), replacing $t$ in (\ref{4.4}) by $h$, we get
\begin{equation}\label{4.8}
\tilde{S}^1(P,q,h)=\sum_{\alpha}\sum_{k=0,1,\cdots}G_{\alpha}^{[k]}(P,q)h^k J_{\alpha}^h=\sum_{\alpha}\sum_{k=0,1,\cdots}G_{\alpha}^{[k]}(P,q)\sum_{\beta\in \Lambda 0_k,\alpha}k!J_{\beta}^h,
\end{equation}
where $0_k$ denotes the index containing $k$ zeros $\underbrace{(0,\cdots,0)}_{k}$. Since $h^k=k!J_{0_k}$, the second equality in (\ref{4.8}) is due to the relation (\cite{kloeden})
\begin{equation}\label{4.9}
\prod_{i=1}^n J_{\alpha_i}=\sum_{\beta\in\Lambda \alpha_1,\cdots,\alpha_n}J_{\beta}.
\end{equation}
Now rearrange the summation terms in (\ref{4.8}) according to different $\beta$, we have
\begin{eqnarray}\label{4.10}
\tilde{S}^1(P,q,h)&=&\sum_{\beta}\left(\sum_{\small \begin{array}{c}k=0,\cdots,l(\beta)-1,\\ \beta\in\Lambda 0_k,\alpha\end{array}}k!\tau_{\beta}(0_k,\alpha)G_{\alpha}^{[k]}(P,q)\right)J_{\beta}^h,\nonumber\\
&=:&\sum_{\beta}\bar{G}_{\beta}(P,q)J_{\beta}^h,
\end{eqnarray}
where $\tau_{\beta}(0_k,\alpha)$ denotes the number of $\beta$ appearing in $\Lambda_{0_k, \alpha}$, and
\begin{equation}\label{4.105}
\bar{G}_{\beta}(P,q):=\sum_{\small \begin{array}{c}k=0,\cdots,l(\beta)-1,\\ \beta\in\Lambda 0_k,\alpha\end{array}}k!\tau_{\beta}(0_k,\alpha)G_{\alpha}^{[k]}(P,q).
\end{equation}

On the other hand, we have the generating function $S^1(P,q,h)=\sum_{\alpha\in \Lambda_{\psi}}G_{\alpha}(P,q)J_{\alpha}^h$ for the numerical method $\psi_h:(p,q)\mapsto(P,Q)$ defined by the relation
\begin{equation*}
 P=p-\frac{\partial S^1}{\partial q} ,\quad Q=q+\frac{\partial S^1}{\partial P},
\end{equation*}
which is equivalent to
\begin{equation*}
\left(\begin{array}{c}P\\Q\end{array}\right)=\left(\begin{array}{c}p\\q\end{array}\right)+J^{-1}\nabla S^1,
\end{equation*}
with $J=\left(\begin{array}{cc}0&1\\-1&0\end{array}\right)$. Thus,
\begin{equation}\label{4.11}
\phi(P,Q)=\phi(p,q)+\nabla\phi(p,q)\cdot J^{-1}\nabla S^1+\frac{1}{2!}\nabla^2\phi(p,q)(J^{-1}\nabla S^1,J^{-1}\nabla S^1)+\cdots .
\end{equation}
Meanwhile, for the modified equation,
\begin{equation*}
\left(\begin{array}{c}\tilde{P}\\\tilde{Q}\end{array}\right)=\left(\begin{array}{c}p\\q\end{array}\right)+J^{-1}\nabla \tilde{S}^1,
\end{equation*}
which implies that
\begin{equation}\label{4.12}
\phi(\tilde{P},\tilde{Q})=\phi(p,q)+\nabla\phi(p,q)\cdot J^{-1}\nabla \tilde{S}^1+\frac{1}{2!}\nabla^2\phi(p,q)(J^{-1}\nabla \tilde{S}^1,J^{-1}\nabla \tilde{S}^1)+\cdots .
\end{equation}
To obtain a modified equation of weak $k$-th order apart from the numerical method $\psi_h$, i.e.,
\begin{equation*}
 | \mathbf{E}(\phi(P,Q))-\mathbf{E}(\phi(\tilde{P},\tilde{Q}))|=O(h^{k+1}),
\end{equation*}
we need to equate the expectation of the right hand side of (\ref{4.11}) and (\ref{4.12}) up to the term of $h^{k}$. Note that, after taking expectations, the terms of order $h^{\frac{1}{2}}$, $h^{\frac{3}{2}}$, etc., vanish. Therefore, for $k=1$, for example, we need to have the coefficients of $h$ to be equal within the following pairs respectively:
\begin{eqnarray}\label{4.13}
&&\mathbf{E}\frac{\partial S^1}{\partial q}\,\,\,\, \mbox{and} \,\,\,\, \mathbf{E}\frac{\partial \tilde{S}^1}{\partial q};\quad \mathbf{E}\frac{\partial S^1}{\partial P}\,\,\,\, \mbox{and} \,\,\,\, \mathbf{E}\frac{\partial \tilde{S}^1}{\partial P};\quad \mathbf{E}\left(\frac{\partial S^1}{\partial q}\right)^2\,\,\,\, \mbox{and} \,\,\,\,\mathbf{E}\left(\frac{\partial \tilde{S}^1}{\partial q}\right)^2;\nonumber \\
 &&\mathbf{E}\frac{\partial S^1}{\partial q}\frac{\partial S^1}{\partial P}\,\,\,\, \mbox{and} \,\,\,\, \mathbf{E}\frac{\partial \tilde{S}^1}{\partial q}\frac{\partial \tilde{S}^1}{\partial P};\quad \mathbf{E}\left(\frac{\partial S^1}{\partial P}\right)^2\,\,\,\, \mbox{and} \,\,\,\, \mathbf{E}\left(\frac{\partial \tilde{S}^1}{\partial P}\right)^2.
\end{eqnarray}
For $k=2$, in addition to equating the coefficients of $h$, as well as those of $h^2$ in the pairs given in (\ref{4.13}), respectively, we need also to have the coefficients of $h^2$ equal in the followig pairs:
\begin{eqnarray}\label{4.14}
&&\mathbf{E}\left(\frac{\partial S^1}{\partial q}\right)^3\,\,\,\, \mbox{and} \,\,\,\, \mathbf{E}\left(\frac{\partial \tilde{S}^1}{\partial q}\right)^3;\quad
\mathbf{E}\left(\frac{\partial S^1}{\partial q}\right)^2\frac{\partial S^1}{\partial P}\,\,\,\, \mbox{and} \,\,\,\, \mathbf{E}\left(\frac{\partial \tilde{S}^1}{\partial q}\right)^2\frac{\partial \tilde{S}^1}{\partial P};\nonumber\\
&&\mathbf{E}\frac{\partial S^1}{\partial q}\left(\frac{\partial S^1}{\partial P}\right)^2\,\,\,\, \mbox{and} \,\,\,\, \mathbf{E}\frac{\partial \tilde{S}^1}{\partial q}\left(\frac{\partial \tilde{S}^1}{\partial P}\right)^2;\quad
 \mathbf{E}\left(\frac{\partial S^1}{\partial P}\right)^3\,\,\,\, \mbox{and} \,\,\,\, \mathbf{E}\left(\frac{\partial \tilde{S}^1}{\partial P}\right)^3;\nonumber\\
 &&\mathbf{E}\left(\frac{\partial S^1}{\partial q}\right)^4\,\,\,\, \mbox{and} \,\,\,\, \mathbf{E}\left(\frac{\partial \tilde{S}^1}{\partial q}\right)^4;\quad
\mathbf{E}\left(\frac{\partial S^1}{\partial q}\right)^3\frac{\partial S^1}{\partial P}\,\,\,\, \mbox{and} \,\,\,\, \mathbf{E}\left(\frac{\partial \tilde{S}^1}{\partial q}\right)^3\frac{\partial \tilde{S}^1}{\partial P};\nonumber\\
&&\mathbf{E}\left(\frac{\partial S^1}{\partial q}\right)^2\left(\frac{\partial S^1}{\partial P}\right)^2\,\,\,\, \mbox{and} \,\,\,\, \mathbf{E}\left(\frac{\partial \tilde{S}^1}{\partial q}\right)^2\left(\frac{\partial \tilde{S}^1}{\partial P}\right)^2;\quad
 \mathbf{E}\frac{\partial S^1}{\partial q}\left(\frac{\partial S^1}{\partial P}\right)^3\,\,\,\, \mbox{and} \nonumber\\&& \mathbf{E}\frac{\partial \tilde{S}^1}{\partial q}\left(\frac{\partial \tilde{S}^1}{\partial P}\right)^3;\,\,\,\,
 \mathbf{E}\left(\frac{\partial S^1}{\partial P}\right)^4\,\,\,\, \mbox{and} \,\,\,\, \mathbf{E}\left(\frac{\partial \tilde{S}^1}{\partial P}\right)^4.
\end{eqnarray}
In general, for $h^k$, we need to equate its coefficients in the expectation of the pairs up to of $2k$-th power of the partial derivatives of $S^1$ and $\tilde{S}^1$. We can then in this process determine the unknown functions $H_{j}^{[i]}$ for $j=0,1,2,\cdots$ and $i=1,2,\cdots$.

The concrete process will be illustrated in the next section.\hfill $\Box$
\section{Some examples}
{\bf Example 1}. A linear stochastic oscillator.

The linear stochastic oscillator (\cite{melbo})
\begin{equation}\label{5.1}\begin{array}{l}
dp=-qdt+\sigma dW(t),\quad p(0)=0,\\
dq=pdt,\quad q(0)=1
\end{array}
\end{equation}
is a stochastic Hamiltonian system with $H=\frac{1}{2}(p^2+q^2)$, $H_1=-\sigma q$. According to (\ref{3.9}) and (\ref{3.9.5}), the coefficient functions $G_{\alpha}$ of the generating function $S^1$ associated with this system are
\begin{equation}\label{5.2}\begin{array}{l}
G_{0}=H=\frac{1}{2}(P^2+q^2),\quad G_1=H_1=-\sigma q,\quad G_{(1,1)}=0,\\
G_{(0,1)}=-\sigma P,\quad G_{(1,0)}=0, \quad G_{(1,1,1)}=0,\quad G_{(0,0)}=Pq,\,\,\cdots.\end{array}
\end{equation}

To obtain a symplectic method of weak order 1, however, we only need to include in the series of $S^1$ those terms of $I_{\alpha}$ with $l(\alpha)\leq 1$. For this, we should use the relations
\begin{equation}\label{5.3}
J_0=I_0,\quad J_1=I_1,\quad J_{(1,1)}=I_{(1,1)}+\frac{1}{2}I_{0}
\end{equation}
to get the correct coefficient of $I_{0}$, which is then $G_0+\frac{1}{2}G_{(1,1)}$, and that of $I_1$ which is $G_1$. Thus, the generating function for a first order weak symplectic scheme is
\begin{equation}\label{5.4}
S^1(P,q,h)=G_1 I_1+(G_0+\frac{1}{2}G_{(1,1)})I_0.
\end{equation}
According to the relation (\ref{3.4}), the weak first order symplectic scheme generated by it is
\begin{equation}\label{5.5}\begin{array}{l}
p_{n+1}=p_n-hq_n+\sigma \Delta W_n,\\
q_{n+1}=q_n+hp_{n+1},\end{array}
\end{equation}
with $\Delta W_n=W(t_{n+1})-W(t_n)$, which is just the stochastic symplectic Euler method. We next find the $\tilde{S}^1=\sum_{\alpha}\bar{G}_{\alpha}J_{\alpha}^h$ associated with this method. According to the formulae (\ref{4.105}), (\ref{4.7}) and (\ref{4.7.5}), we have
\begin{equation}\label{5.6}\begin{array}{l}
\bar{G}_0=G_0^{[0]}=H_0^{[0]}=H=\frac{1}{2}(P^2+q^2),\quad \bar{G}_1=G_1^{[0]}=H_1^{[0]}=H_1=-\sigma q,\\
\bar{G}_{(1,1)}=G_{(1,1)}^{[0]}=\frac{\partial H_1^{[0]}}{\partial q}\frac{\partial G_1^{[0]}}{\partial p}=0,\\
\bar{G}_{(0,1)}=G_{(0,1)}^{[0]}+G_1^{[1]}=\frac{\partial H_1^{[0]}}{\partial q}\frac{\partial G_0^{[0]}}{\partial p}+H_1^{[1]}=-\sigma P+H_1^{[1]},\\
\bar{G}_{(1,0)}=G_{(1,0)}^{[0]}+G_1^{[1]}=\frac{\partial H_0^{[0]}}{\partial q}\frac{\partial G_1^{[0]}}{\partial p}+H_1^{[1]}=H_1^{[1]},\\
\bar{G}_{(0,0)}=G_{(0,0)}^{[0]}+2G_0^{[1]}=\frac{\partial H_0^{[0]}}{\partial q}\frac{\partial G_0^{[0]}}{\partial p}+2H_0^{[1]}=Pq+2H_0^{[1]},\\
\bar{G}_{(1,1,0)}=G_{(1,1,0)}^{[0]}+G_{(1,1)}^{[1]}\\ \qquad \quad=\frac{\partial H_0^{[0]}}{\partial q}\frac{\partial G_{(1,1)}^{[0]}}{\partial p}+\frac{\partial^2 H_0^{[0]}}{\partial q^2}\frac{\partial G_1^{[0]}}{\partial p}^2+\frac{\partial H_1^{[0]}}{\partial q}\frac{\partial G_1^{[1]}}{\partial p}+\frac{\partial H_1^{[1]}}{\partial q}\frac{\partial G_1^{[0]}}{\partial p}\\ \qquad \quad =-\sigma \frac{\partial H_1^{[1]}}{\partial p}\\
\bar{G}_{(0,1,1)}=G_{(0,1,1)}^{[0]}+G_{(1,1)}^{[1]}\\ \qquad \quad=\frac{\partial H_1^{[0]}}{\partial q}\frac{\partial G_{(0,1)}^{[0]}}{\partial p}+\frac{\partial^2 H_1^{[0]}}{\partial q^2}\frac{\partial G_1^{[0]}}{\partial p}\frac{\partial G_0^{[0]}}{\partial p}+\frac{\partial H_1^{[0]}}{\partial q}\frac{\partial G_1^{[1]}}{\partial p}+\frac{\partial H_1^{[1]}}{\partial q}\frac{\partial G_1^{[0]}}{\partial p}\\ \qquad \quad =\sigma^2-\sigma \frac{\partial H_1^{[1]}}{\partial p},\\
\bar{G}_{(1,1,1)}=G_{(1,1,1)}^{[0]}=\frac{\partial H_1^{[0]}}{\partial q}\frac{\partial G_{(1,1)}^{[0]}}{\partial p}+\frac{\partial^2 H_1^{[0]}}{\partial q^2}\frac{\partial G_1^{[0]}}{\partial p}^2=0,\\
\bar{G}_{(1,1,1,1)}=G_{(1,1,1,1)}^{[0]}=\frac{\partial H_1^{[0]}}{\partial q}\frac{\partial G_{(1,1,1)}^{[0]}}{\partial p}+3\frac{\partial^2 H_1^{[0]}}{\partial q^2}\frac{\partial G_{(1,1)}^{[0]}}{\partial p}\frac{\partial G_1^{[0]}}{\partial p}+\frac{\partial^3 H_1^{[0]}}{\partial q^3}\frac{\partial G_1^{[0]}}{\partial p}^3=0,\\
\vdots
\end{array}
\end{equation}

It is easy to check that the coefficients of $h$ within each pair in (\ref{4.13}) are naturally equal, since $G_0=\bar{G}_0$, $G_1=\bar{G}_1$, and $G_{(1,1)}=\bar{G}_{(1,1)}=0$. Equating the coefficients of $h^2$ within each pair in (\ref{4.13}) and (\ref{4.14}) gives
\begin{eqnarray}
&&\frac{1}{2}\sigma\frac{\partial^2H_1^{[1]}}{\partial P\partial q}-\frac{\partial H_0^{[1]}}{\partial q}=\frac{1}{2}P,\label{5.7.1}\\
&&\frac{1}{2}\sigma\frac{\partial^2H_1^{[1]}}{\partial P^2}-\frac{\partial H_0^{[1]}}{\partial P}=\frac{1}{2}q,\label{5.7.2}\\
&&\frac{\partial H_1^{[1]}}{\partial q}=0,\label{5.7.3}\\
&&\frac{\partial H_1^{[1]}}{\partial P}=\frac{1}{2}\sigma.\label{5.7.3}
\end{eqnarray}
Substituting (\ref{5.7.3}) into (\ref{5.7.1}) and (\ref{5.7.2}) gives
\begin{equation}\label{5.11}
\frac{\partial H_0^{[1]}}{\partial P}=-\frac{1}{2}q,\quad \frac{\partial H_0^{[1]}}{\partial q}=-\frac{1}{2}P.
\end{equation}
Thus, according to (\ref{4.1})-(\ref{4.2}), the modified equation of weak second order apart from the numerical method (\ref{5.5}) is
\begin{eqnarray}\label{5.12}
dp&=&(-q+h\frac{p}{2})dt+\sigma dW(t),\nonumber\\
dq&=&(p-h\frac{q}{2})dt+h\frac{\sigma}{2} dW(t),
\end{eqnarray}
which coincides with the result (Eq. (4.16) in \cite{zy}) about the modified equation of the symplectic Euler method for a Langevin equation (Eq. (4.14) in \cite{zy}) as $V'(q)=q$ and $\gamma=0$.

For modified equation of weak third order apart from the numerical method (\ref{5.5}), we need to equate the coefficients of $h^3$ within corresponding pairs to determine more unknown coefficients $H_j^{[i]}$, and so on and so forth for even higher orders.
\\ \\
{\bf Example 2}. A model for synchrotron oscillations of particles in storage rings

The model (\cite{mil2}) is
\begin{eqnarray}\label{5.13}
dp&=&-\omega^2 \sin q dt-\sigma_1 \cos q \circ dW_1(t)-\sigma_2 \sin q\circ dW_2(t),\nonumber\\
dq&=&pdt,
\end{eqnarray}
where $p$ and $q$ are scaler. It is a stochastic Hamiltonian system with
\begin{equation}\label{5.14}
H_0=-\omega^2\cos q+\frac{1}{2}p^2,\quad H_1=\sigma_1\sin q,\quad H_2=-\sigma_2\cos q.
\end{equation}
According to (\ref{3.9}) and (\ref{3.9.5}),
\begin{eqnarray}\label{5.15}
G_0&=&H_0=-\omega^2\cos q+\frac{1}{2}p^2,\quad G_1=H_1=\sigma_1\sin q,\nonumber\\
G_2&=&H_2=-\sigma_2\cos q,\quad G_{(1,1)}=G_{(2,2)}=0,\cdots.
\end{eqnarray}
Thus the generating function $S^1$ for a weak first order symplectic scheme is
\begin{equation}\label{5.16}
S^1(P,q,h)=(G_0+\frac{1}{2}G_{(1,1)}+\frac{1}{2}G_{(2,2)} )I_0+G_{1}I_1+G_2I_2,
\end{equation}
and the scheme generated by it via the relation (\ref{3.4}) is
\begin{eqnarray}\label{5.17}
p_{n+1}&=&p_n-(h\omega^2\sin q_n+\Delta_n W_1 \sigma_1 \cos q_n+\Delta_n W_2 \sigma_2 \sin q_n),\nonumber\\
q_{n+1}&=& q_n+hp_{n+1},
\end{eqnarray}
where $\Delta_nW_i=W_i(t_n+1)-W_i(t_n)$ for $i=1,2$. For the generating function $\tilde{S}^1(P,q,t)$ of the modified equation of (\ref{5.17}) at $t=h$ we have $\tilde{S}^1(P,q,h)=\sum_{\alpha}\bar{G}_{\alpha}(P,q)J_{\alpha}^h$, where, according to the formulae (\ref{4.105}), (\ref{4.7}) and (\ref{4.7.5}),
\begin{equation}\label{5.18}\begin{array}{l}
\bar{G}_0=G_0^{[0]}=H_0^{[0]}=H_0=-\omega^2\cos q+\frac{1}{2}P^2,\quad \bar{G}_1=G_1^{[0]}=H_1^{[0]}=H_1=\sigma_1\sin q,\\
\bar{G}_2=G_1^{[0]}=H_2^{[0]}=H_2=-\sigma_2\cos q,\quad \bar{G}_{(1,1)}=G_{(1,1)}^{[0]}=\frac{\partial H_1^{[0]}}{\partial q}\frac{\partial G_1^{[0]}}{\partial p}=0,\\
\bar{G}_{(1,2)}=G_{(1,2)}^{[0]}=\frac{\partial H_2^{[0]}}{\partial q}\frac{\partial G_1^{[0]}}{\partial p}=0,\quad \bar{G}_{(2,1)}=G_{(2,1)}^{[0]}=\frac{\partial H_1^{[0]}}{\partial q}\frac{\partial G_2^{[0]}}{\partial p}=0,\\
\bar{G}_{(2,2)}=G_{(2,2)}^{[0]}=\frac{\partial H_2^{[0]}}{\partial q}\frac{\partial G_2^{[0]}}{\partial p}=0,\\
\bar{G}_{(0,1)}=G_{(0,1)}^{[0]}+G_1^{[1]}=\frac{\partial H_1^{[0]}}{\partial q}\frac{\partial G_0^{[0]}}{\partial p}+H_1^{[1]}=\sigma_1 P\cos q+H_1^{[1]},\\
\bar{G}_{(0,2)}=G_{(0,2)}^{[0]}+G_2^{[1]}=\frac{\partial H_2^{[0]}}{\partial q}\frac{\partial G_0^{[0]}}{\partial p}+H_2^{[1]}=\sigma_2 P\sin q+H_2^{[1]},\\
\bar{G}_{(1,0)}=G_{(1,0)}^{[0]}+G_1^{[1]}=\frac{\partial H_0^{[0]}}{\partial q}\frac{\partial G_1^{[0]}}{\partial p}+H_1^{[1]}=H_1^{[1]},\\
\bar{G}_{(2,0)}=G_{(2,0)}^{[0]}+G_2^{[1]}=\frac{\partial H_0^{[0]}}{\partial q}\frac{\partial G_2^{[0]}}{\partial p}+H_2^{[1]}=H_2^{[1]},\\
\bar{G}_{(1,1,1)}=G_{(1,1,1)}^{[0]}=\frac{\partial H_1^{[0]}}{\partial q}\frac{\partial G_{(1,1)}^{[0]}}{\partial p}+\frac{\partial^2H_1^{[0]}}{\partial q^2}\left(\frac{\partial G_1^{[0]}}{\partial P}\right)^2=0,\\
\bar{G}_{(1,1,2)}=G_{(1,1,2)}^{[0]}=\frac{\partial H_2^{[0]}}{\partial q}\frac{\partial G_{(1,1)}^{[0]}}{\partial p}+\frac{\partial^2H_2^{[0]}}{\partial q^2}\left(\frac{\partial G_1^{[0]}}{\partial P}\right)^2=0,\\
\bar{G}_{(1,2,1)}=G_{(1,2,1)}^{[0]}=\frac{\partial H_1^{[0]}}{\partial q}\frac{\partial G_{(1,2)}^{[0]}}{\partial p}+\frac{\partial^2H_1^{[0]}}{\partial q^2}\frac{\partial G_1^{[0]}}{\partial P}\frac{\partial G_2^{[0]}}{\partial P}=0,\\
\bar{G}_{(1,2,2)}=G_{(1,2,2)}^{[0]}=\frac{\partial H_2^{[0]}}{\partial q}\frac{\partial G_{(1,2)}^{[0]}}{\partial p}+\frac{\partial^2H_2^{[0]}}{\partial q^2}\frac{\partial G_1^{[0]}}{\partial P}\frac{\partial G_2^{[0]}}{\partial P}=0,\\
\bar{G}_{(2,1,1)}=G_{(2,1,1)}^{[0]}=\frac{\partial H_1^{[0]}}{\partial q}\frac{\partial G_{(2,1)}^{[0]}}{\partial p}+\frac{\partial^2H_1^{[0]}}{\partial q^2}\frac{\partial G_2^{[0]}}{\partial P}\frac{\partial G_1^{[0]}}{\partial P}=0,\\
\bar{G}_{(2,1,2)}=G_{(2,1,2)}^{[0]}=\frac{\partial H_2^{[0]}}{\partial q}\frac{\partial G_{(2,1)}^{[0]}}{\partial p}+\frac{\partial^2H_2^{[0]}}{\partial q^2}\frac{\partial G_2^{[0]}}{\partial P}\frac{\partial G_1^{[0]}}{\partial P}=0,\\
\bar{G}_{(2,2,1)}=G_{(2,2,1)}^{[0]}=\frac{\partial H_1^{[0]}}{\partial q}\frac{\partial G_{(2,2)}^{[0]}}{\partial p}+\frac{\partial^2H_1^{[0]}}{\partial q^2}\left(\frac{\partial G_2^{[0]}}{\partial P}\right)^2=0,\\
\bar{G}_{(2,2,2)}=G_{(2,2,2)}^{[0]}=\frac{\partial H_2^{[0]}}{\partial q}\frac{\partial G_{(2,2)}^{[0]}}{\partial p}+\frac{\partial^2H_2^{[0]}}{\partial q^2}\left(\frac{\partial G_2^{[0]}}{\partial P}\right)^2=0,\\
\bar{G}_{(0,0)}=G_{(0,0)}^{[0]}+2G_0^{[1]}=\frac{\partial H_0^{[0]}}{\partial q}\frac{\partial G_0^{[0]}}{\partial p}+2H_0^{[1]}=\omega^2 P\sin q+2H_0^{[1]},\\
\bar{G}_{(1,1,0)}=G_{(1,1,0)}^{[0]}+G_{(1,1)}^{[1]}=\sigma_1\cos q \frac{\partial H_1^{[1]}}{\partial p},
\\ \bar{G}_{(1,2,0)}=G_{(1,2,0)}^{[0]}+G_{(1,2)}^{[1]}=\sigma_2\sin q \frac{\partial H_1^{[1]}}{\partial P}, \\
\bar{G}_{(2,1,0)}=G_{(2,1,0)}^{[0]}+G_{(2,1)}^{[1]}=\sigma_1\cos q \frac{\partial H_2^{[1]}}{\partial p}, \\
\bar{G}_{(2,2,0)}=G_{(2,2,0)}^{[0]}+G_{(2,2)}^{[1]}=\sigma_2\sin q \frac{\partial H_2^{[1]}}{\partial P}, \\
\bar{G}_{(0,1,1)}=G_{(0,1,1)}^{[0]}+G_{(1,1)}^{[1]}=\sigma_1^2\cos^2 q+\sigma_1\cos q \frac{\partial H_1^{[1]}}{\partial p},\\ \bar{G}_{(0,1,2)}=G_{(0,1,2)}^{[0]}+G_{(1,2)}^{[1]}=\sigma_1\sigma_2\sin q\cos q+\sigma_2\sin q \frac{\partial H_1^{[1]}}{\partial p},
\end{array}
\end{equation}
\begin{equation}\label{5.19}
\begin{array}{l}
\bar{G}_{(0,2,1)}=G_{(0,2,1)}^{[0]}+G_{(2,1)}^{[1]}=\sigma_1\sigma_2\sin q\cos q+\sigma_1\cos q \frac{\partial H_2^{[1]}}{\partial p},\\ \bar{G}_{(0,2,2)}=G_{(0,2,2)}^{[0]}+G_{(2,2)}^{[1]}=\sigma_2^2\sin^2 q\cos q+\sigma_2\sin q \frac{\partial H_2^{[1]}}{\partial p},\\
\bar{G}_{(1,0,1)}=G_{(1,0,1)}^{[0]}+G_{(1,1)}^{[1]}=\sigma_1\cos q\frac{\partial H_1^{[1]}}{\partial P},\\
\bar{G}_{(1,0,2)}=G_{(1,0,2)}^{[0]}+G_{(1,2)}^{[1]}=\sigma_2\sin q\frac{\partial H_1^{[1]}}{\partial P},\\
\bar{G}_{(2,0,1)}=G_{(2,0,1)}^{[0]}+G_{(2,1)}^{[1]}=\sigma_1\cos q\frac{\partial H_2^{[1]}}{\partial P},\\
\bar{G}_{(2,0,2)}=G_{(2,0,2)}^{[0]}+G_{(2,2)}^{[1]}=\sigma_2\sin q\frac{\partial H_2^{[1]}}{\partial P},\\
\bar{G}_{(1,1,1,1)}=G_{(1,1,1,1)}^{[0]}=\frac{\partial H_1^{[0]}}{\partial q}\frac{\partial G_{(1,1,1)}^{[0]}}{\partial p}+3\frac{\partial^2 H_1^{[0]}}{\partial q^2}\frac{\partial G_{(1,1)}^{[0]}}{\partial p}\frac{\partial G_1^{[0]}}{\partial p}+\frac{\partial^3 H_1^{[0]}}{\partial q^3}\left(\frac{\partial G_1^{[0]}}{\partial p}\right)^3=0,\\
\bar{G}_{(1,1,1,2)}=G_{(1,1,1,2)}^{[0]}=0,\quad \bar{G}_{(1,1,2,1)}=G_{(1,1,2,1)}^{[0]}=0,\\
\bar{G}_{(1,1,2,2)}=G_{(1,1,2,2)}^{[0]}=0,\quad \bar{G}_{(1,2,1,1)}=G_{(1,2,1,1)}^{[0]}=0,\\
\bar{G}_{(1,2,1,2)}=G_{(1,2,1,2)}^{[0]}=0,\quad \bar{G}_{(1,2,2,1)}=G_{(1,2,2,1)}^{[0]}=0,\\
\bar{G}_{(1,2,2,2)}=G_{(1,2,2,2)}^{[0]}=0,\quad \bar{G}_{(2,1,1,1)}=G_{(2,1,1,1)}^{[0]}=0,\\
\bar{G}_{(2,1,1,2)}=G_{(2,1,1,2)}^{[0]}=0,\quad \bar{G}_{(2,1,2,1)}=G_{(2,1,2,1)}^{[0]}=0,\\
\bar{G}_{(2,1,2,2)}=G_{(2,1,2,2)}^{[0]}=0,\quad \bar{G}_{(2,2,1,1)}=G_{(2,2,1,1)}^{[0]}=0,\\
\bar{G}_{(2,2,1,2)}=G_{(2,2,1,2)}^{[0]}=0,\quad \bar{G}_{(2,2,2,1)}=G_{(2,2,2,1)}^{[0]}=0,\\
\bar{G}_{(2,2,2,2)}=G_{(2,2,2,2)}^{[0]}=0,\\
\vdots
\end{array}
\end{equation}
Due to $G_0=\bar{G}_0$, $G_1=\bar{G}_1$, $G_2=\bar{G}_2$, $\bar{G}_{(1,1)}=\bar{G}_{(2,2)}=0$, the coefficients of $h$ in each pair of (\ref{4.13}) are equal. Equating the coefficients of $h^2$ in each pair of (\ref{4.13}) and (\ref{4.14}), we obtain
\begin{eqnarray}
&&\frac{\partial H_0^{[1]}}{\partial q}+\frac{1}{2}\left(-\sigma_1\sin q\frac{\partial H_1^{[1]}}{\partial P}+\sigma_1\cos q\frac{\partial^2 H_1^{[1]}}{\partial P\partial q}\right)\nonumber\\
&+&\frac{1}{2}\left(\sigma_2\cos q\frac{\partial H_2^{[1]}}{\partial P}+\sigma_2\sin q \frac{\partial^2H_2^{[1]}}{\partial P\partial q}\right)\nonumber\\
&=&-\frac{1}{2}\omega^2 P\cos q+\frac{1}{2}\sigma_1^2\sin q\cos q-\frac{1}{2}\sigma_2^2\sin q\cos q,\label{5.20}\\
&&\frac{\partial H_0^{[1]}}{\partial P}+\frac{1}{2}\sigma_1\cos q\frac{\partial^2 H_1^{[1]}}{\partial P^2}+\frac{1}{2}\sigma_2\sin q\frac{\partial^2 H_2^{[1]}}{\partial P^2}
=-\frac{1}{2}\omega^2 \sin q,\label{5.21}\\
&&2\left(\sigma_1\cos q\frac{\partial H_1^{[1]}}{\partial q}+\sigma_2\sin q\frac{\partial H_2^{[1]}}{\partial q}\right)=(\sigma_1^2-\sigma_2^2)P\sin q\cos q,\label{5.22}\\
&&2\left(\sigma_1\cos q\frac{\partial H_1^{[1]}}{\partial P}+\sigma_2\sin q\frac{\partial H_2^{[1]}}{\partial P}\right)=-\sigma_1^2\cos^2 q-\sigma_2^2\sin^2 q.\label{5.23}
\end{eqnarray}
Form (\ref{5.22}) and (\ref{5.23}) it follows that
\begin{eqnarray}\label{5.24}
\frac{\partial H_1^{[1]}}{\partial P}&=&-\frac{1}{2}\sigma_1\cos q,\quad \frac{\partial H_1^{[1]}}{\partial q}=\frac{1}{2}\sigma_1P\sin q,\nonumber\\
\frac{\partial H_2^{[1]}}{\partial P}&=&-\frac{1}{2}\sigma_2\sin q,\quad \frac{\partial H_2^{[1]}}{\partial q}=-\frac{1}{2}\sigma_2P\cos q.
\end{eqnarray}
Substituting (\ref{5.24}) into (\ref{5.20}) and (\ref{5.21}), we get
\begin{equation}\label{5.25}
\frac{\partial H_0^{[1]}}{\partial q}=-\frac{1}{2}\omega^2 P\cos q,\quad \frac{\partial H_0^{[1]}}{\partial P}=-\frac{1}{2}\omega^2 \sin q.
\end{equation}
Thus, the modified equation of weak second order apart from the numerical method (\ref{5.17}) is
\begin{eqnarray}\label{5.26}
dp&=&(-\omega^2\sin q+\frac{h}{2}\omega^2 p\cos q)dt-(\sigma_1\cos q+\frac{h}{2}\sigma_1 p\sin q)\circ dW_1(t)\nonumber\\&-&(\sigma_2\sin q-\frac{h}{2}\sigma_2 p\cos q)\circ dW_2(t),\nonumber\\
dq&=&(p-\frac{h}{2}\omega^2\sin q)dt-\frac{h}{2}\sigma_1\cos q\circ dW_1(t)-\frac{h}{2}\sigma_2\sin q\circ dW_2(t).
\end{eqnarray}
For higher order modified equations, more $\bar{G}_{\alpha}$ should be computed, and more equations should be satisfied, which increases the computational complexity.
\\ \\
{\bf Remark}. In the aforementioned examples, the Hamiltonian functions associated with the Gaussian noises are `half' multiplicative, which means that they either depend only on $q$ or only on $p$, such as $H_1=\sigma_1\sin q$, $H_2=-\sigma_2\cos q$ in Example 2. In the case of `totally' multiplicative Hamiltonian functions, i.e., $H_i\,\,(i>0)$ depend on both $p$ and $q$, it is not possible to write out a modified equation of second order accuracy in the form of (\ref{4.1})-(\ref{4.2}) using the procedure given above. An example for this is the Kubo oscillator (\cite{mil2})
\begin{eqnarray}\label{kubo}
dp&=&-aqdt-\sigma q\circ dW(t),\quad p(0)=p_0,\nonumber\\
dq&=&apdt+\sigma p\circ dW(t),\quad q(0)=q_0,
\end{eqnarray}
where $H_1=\frac{\sigma}{2}(p^2+q^2)$. This coincides with the assertion in \cite{sh} that it is impossible to develop modified equations of higher order accuracy ($\geq 2$) for SDEs with multiplicative noises by working only with deterministic perturbations of the drift and diffusion coefficients. In such cases we may need to assume the perturbed Hamiltonian functions $\tilde{H}_i\,\,(i\geq 0)$ in (\ref{4.1}) to have a more general and delicate formulation.
\section{Conclusion}
Stochastic symplectic methods for stochastic Hamiltonian systems are shown to have good long time numerical behavior (\cite{milbook,mil1,mil2}). The reason for this is expected to be revealed by stochastic backward error analysis, as in the deterministic cases. In this paper, we prove that the modified equations of the weak stochastic symplectic methods for stochastic Hamiltonian systems (SHS) with additive or `half' multiplicative noises are perturbed SHS of the original SHS, which gives insight to their superiority in numerical simulations. Meanwhile, we propose the procedure of constructing the modified equations up to certain weak order of accuracy for weak stochastic symplectic methods using their generating functions, instead of the existing method of writing stochastic modified equations. Applications of the method to two examples succeed in rewriting the modified equation obtained by other methods in existing literature, and in establishing a second weak order modified equation of a stochastic symplectic method for a `half' multiplicative SHS. These give validation to the proposed theory and method.

\section*{Acknowledgements}
The first author is supported by the NNSFC (No. 11071251, No. 11471310), and by the 2013 Director Foundation of UCAS. The second author is supported by the NNSFC (No. 91130003, No. 11021101 and No. 11290142).



\end{document}